\documentstyle[12pt,amssym]{article}

\font\elevenof=msbm10 at 11pt
\font\fourteenof=msbm10 at 14pt
\font\twentyof=msbm10 at 20pt

\def\R{\mbox{$\Bbb R$}}
\def\N{\mbox{$\Bbb N$}}
\def\case#1#2{{\textstyle{#1\over #2}}}
\def\zb{\overline{z}}
\def\Tb{\overline{T}}
\def\cqfd{\hfill$\blacksquare$}

\newcommand{\bin}[2]{\left[ \begin{array}{c} {#1} \\ {#2} \end{array} \right]_q}
\newcommand{\Res}{\mathop{\rm Res}\nolimits}

\renewcommand{\theequation}{\arabic{section}.\arabic{equation}}

\renewcommand{\Re}{\mathop{\rm Re}\nolimits}
\renewcommand{\Im}{\mathop{\rm Im}\nolimits}
\newtheorem{proposition}{Proposition}[section]
\newtheorem{corollary}[proposition]{Corollary}
\newtheorem{lemma}[proposition]{Lemma}

\title{
\hfill{\normalsize ULB/229/CQ/98/2}\\
\vspace{1cm}
Unitary representations of the quantum algebra su$_q$(2) on a
real two-dimensional sphere for $q \in \mbox{\twentyof R}^+$ or generic $q \in
S^1$}
\author{M. Irac-Astaud\thanks{E-mail: mici@ccr.jussieu.fr}\\
{\small \sl Laboratoire de Physique Th\'eorique de la Mati\`ere Condens\'ee,
Universit\'e Paris VII,}\\
{\small \sl 2, place Jussieu, F-75251 Paris Cedex 05, France}\\[0.5cm]
C. Quesne\thanks{Directeur de recherches FNRS; E-mail: cquesne@ulb.ac.be} \\
{\small \sl Physique Nucl\'eaire Th\'eorique et Physique Math\'ematique,
Universit\'e Libre de Bruxelles,} \\
{\small \sl Campus de la Plaine CP229, Boulevard~du Triomphe, B-1050 Brussels,
Belgium}}
\date{ }
\begin{document}
\baselineskip=22pt plus 1pt minus 1pt
\maketitle

\begin{abstract}
Some time ago, Rideau and Winternitz introduced a realization of the quantum
algebra~su$_q$(2) on a real two-dimensional sphere, or a real plane, and
constructed a basis for its representations in terms of $q$-special
functions, which
can be expressed in terms of $q$-Vilenkin functions, and are related to little
$q$-Jacobi functions, $q$-spherical functions, and $q$-Legendre polynomials.
In their study, the values of~$q$ were implicitly restricted to $q \in
\mbox{\elevenof R}^+$. In the present paper, we extend their work to the case of
generic values of $q \in S^1$ (i.e., $q$~values different from a root of
unity). In
addition, we unitarize the representations for both types of $q$~values, $q \in
\mbox{\elevenof R}^+$ and generic $q \in S^1$,  by determining some appropriate
scalar products. From the latter, we deduce the orthonormality relations
satisfied
by the $q$-Vilenkin functions.
\end{abstract}

\vspace{0.5cm}

\hspace*{0.3cm}
PACS: 02.30.Gp, 02.20.Sv, 03.65.Fd

\hspace*{0.3cm}
Running title: Unitary representations of quantum algebra
\newpage
%
%
\section{INTRODUCTION} \label{sec:intro}
As is well known, most special functions of mathematical physics admit
extensions to a base~$q$, which are called $q$-special
functions~\cite{exton,andrews,gasper}. In the same way as Lie algebras and their
representations provide a unifying framework for the former, quantum
algebras~\cite{chari} are relevant to the study of the latter (see
e.g.~\cite{vilenkin2} and references quoted therein).\par
%
%
Some time ago, Rideau and Winternitz~\cite{rideau} introduced a realization
of the
quantum algebra~su$_q$(2) on a real sphere~$S^2$ (or, via a stereographic
projection, on a real plane), and constructed a basis for its irreducible
representations (irreps) in terms of some functions $\Psi^J_{MNq}(\theta,\phi)
\propto P^J_{MNq}(\cos\theta) \exp\left(-i(M+N)\phi\right)$. The functions
$P^J_{MNq}(\cos\theta)$ were called $q$-Vilenkin functions because, for $q=1$,
they reduce to functions $P^J_{MN}(\cos\theta)$ introduced by
Vilenkin~\cite{vilenkin1,footnote}, and related to Jacobi polynomials.\par
%
%
Rideau and Winternitz did establish various interesting results for the
$q$-Vilenkin functions, including their recursion relations, explicit
expression,
generating function, and symmetry relations. They also compared them with other
$q$-special functions, such as $q$-hypergeometric series, little $q$-Jacobi
functions, $q$-spherical functions, and $q$-Legendre polynomials. Recently, the
latter polynomials were further studied by Schmidt along similar
lines~\cite{schmidt}.\par
%
%
The realization of~su$_q$(2) on~$S^2$, introduced by Rideau and Winternitz, was
used by one of the present authors~(MIA) to set up su$_q$(2)-invariant
Schr\"odinger equations in the usual framework of quantum
mechanics~\cite{irac1}.
The corresponding radial equations can be easily solved for the ``free''
su$_q$(2)-invariant particle~\cite{irac1}, as well as for the
Coulomb~\cite{irac1}
and oscillator~\cite{irac2} potentials.\par
%
%
Although not explicitly stated in Ref.~\cite{rideau}, the values of the
deformation
parameter~$q$, considered there, are restricted to $q \in \R^+$. Close
examination
indeed shows that the explicit form of the function $Q_{Jq}(\eta)$, $\eta \equiv
\cot^2(\theta/2)$, entering the definition of the $q$-Vilenkin
functions~\cite{rideau}, is not valid for half-integer $J$~values, whenever
$q$ runs
over the unit circle.\par
%
%
Though important both from the $q$-special function viewpoint, and from that of
their applications in quantum mechanics, the question of the su$_q$(2) irrep
unitarity was also left unsolved by Rideau and Winternitz. They only
noticed~\cite{rideau} that their realization of~su$_q$(2) on~$S^2$ is not
unitary
with respect to the scalar product used to unitarize the corresponding
realization
of~su(2), and that a new scalar product should therefore be determined to cope
with this drawback.\par
%
%
The purpose of the present paper is twofold: firstly, to find a solution for
$Q_{Jq}(\eta)$ for generic $q \in S^1$ (i.e., for $q$ different from a root
of unity),
and secondly, to unitarize the representations for both $q\in \R^+$, and
generic $q
\in S^1$. As a consequence, the explicit orthonormality relations of the $q$-Vilenkin
and related functions will be established.\par
%
In Sec.~\ref{sec:representations}, the representations of~su$_q$(2) on~$S^2$,
derived by Rideau and Winternitz, are briefly reviewed. The function
$Q_{Jq}(\eta)$
is determined in Sec.~\ref{sec:Q}. The unitarization of the representations
is dealt
with in Sec.~\ref{sec:unitarization}. Sec.~\ref{sec:conclusion} contains the
conclusion.\par
%
%
\section{REPRESENTATIONS OF su$_q$(2) ON $S^2$} \label{sec:representations}
Let us consider functions $f(\theta,\phi)$ on a sphere~$S^2$, defined by
$x_0^2 +
y_0^2 + z_0^2 = 1/4$. These functions can also be viewed as functions on a real
plane, via the stereographic projection $x = x_0/(1/2-z_0)$, $y =
y_0/(1/2-z_0)$.
In terms of spherical coordinates on~$S^2$ and polar ones on the plane, we have
\begin{eqnarray}
  x_0 & = & \case{1}{2} \sin\theta \cos\phi, \qquad y_0 = \case{1}{2} \sin\theta
         \sin\phi, \qquad z_0 = \case{1}{2} \cos\theta, \nonumber \\
  x & = & \rho \cos\phi, \qquad y = \rho \sin\phi, \qquad \rho =
         \cot\frac{\theta}{2}, \nonumber \\
  0 & \le & \theta \le \pi, \qquad 0 \le \phi < 2\pi, \qquad 0 \le \rho <
\infty.
         \label{eq:coord}
\end{eqnarray}
Instead of the real variables~$x$ and~$y$, one can use complex ones
\begin{equation}
  z = x + iy = \rho\, e^{i\phi}, \qquad \zb = x - iy = \rho\, e^{-i\phi}.
  \label{eq:complex}
\end{equation}
{}Functions $f(\theta,\phi)$ on~$S^2$ can thus be projected onto functions
$f(\rho,\phi)$ on the real plane, or functions $f(z,\zb)$ of a complex
variable and
its conjugate.\par
%
%
The su$_q$(2) generators $H_3$, $H_+$, $H_-$ satisfy the commutation
relations~\cite{chari}
\begin{equation}
  \left[H_3, H_{\pm}\right] = \pm H_{\pm}, \qquad \left[H_+, H_-\right] =
  \left[2H_3\right]_q \equiv \frac{q^{2H_3} - q^{-2H_3}}{q - q^{-1}},
\label{eq:com}
\end{equation}
and the Hermiticity properties
\begin{equation}
  H_3^{\dagger} = H_3, \qquad H_{\pm}^{\dagger} = H_{\mp},  \label{eq:hermite}
\end{equation}
where in Eq.~(\ref{eq:com}), we assume $q = e^{\tau} \in \R^+$, or $q =
e^{i\tau} \in
S^1$ (but different from a root of unity). From~$H_3$ and~$H_{\pm}$, one can
construct a Casimir operator
\begin{equation}
  {\cal C} = H_+ H_- + \left[H_3\right]_q \left[H_3-1\right]_q = H_- H_+
  + \left[H_3\right]_q \left[H_3+1\right]_q,   \label{eq:casimir}
\end{equation}
such that $\left[{\cal C}, H_3\right] = \left[{\cal C}, H_{\pm}\right] = 0$.\par
%
%
The generators $H_3$, $H_+$, $H_-$ can be realized~\cite{rideau} by the
following
operators, acting on functions $f(z,\zb)$ or $f(\theta,\phi)$,
\begin{eqnarray}
  H_3 & = & - z \partial_z + \zb \partial_{\zb} - N = i \partial_{\phi} - N,
          \nonumber \\
  H_+ & = & - z^{-1}\, [T]_q\, q^{\Tb - (N/2)} - q^{T + (N/2)}\, \zb \left[\Tb -
          N\right]_q, \nonumber \\
  H_- & = & z\, [T + N]_q\, q^{\Tb - (N/2)} + q^{T + (N/2)}\, \zb^{-1}
          \left[\Tb\right]_q,  \label{eq:su-q}
\end{eqnarray}
where
\begin{equation}
  T = z \partial_z = - \case{1}{2} \left(\sin\theta \partial_{\theta} + i
  \partial_{\phi}\right), \qquad \Tb = \zb \partial_{\zb} = - \case{1}{2}
  \left(\sin\theta \partial_{\theta} - i \partial_{\phi}\right).
\end{equation}
{}For future use, it is also convenient to write $H_{\pm}$ in terms of polar
coordinates on the real plane as
\begin{equation}
  H_{\pm} = \mp \frac{e^{\mp i\phi}}{q-q^{-1}} \left\{\left(\rho +
\frac{1}{\rho}
  \right) q^{\rho \partial_{\rho} \mp (N/2)} - \rho\, q^{\mp i
\partial_{\phi} \pm
  (3N/2)} - \frac{1}{\rho}\, q^{\pm i \partial_{\phi} \mp (N/2)}\right\}.
\end{equation}
\par
%
%
Basis functions $\Psi^J_{MNq}(z,\zb)$ for the ($2J+1$)-dimensional irrep
of~su$_q$(2) satisfy the relations~\cite{chari}
\begin{eqnarray}
  H_3 \Psi^J_{MNq} & = & M \Psi^J_{MNq}, \qquad H_{\pm} \Psi^J_{MNq} = \left(
         [J\mp M]_q [J\pm M+1]_q\right)^{1/2} \Psi^J_{M\pm1,Nq}, \nonumber \\
  {\cal C} \Psi^J_{MNq} & = & [J]_q [J+1]_q \Psi^J_{MNq},\qquad  M = \{ -J,-J+1,
         \cdots,J\} , \qquad |N| \leq J,   \label{eq:irrep}
\end{eqnarray}
where $J$, $M$ and $N$ are simultaneously integers or half-integers. Let us
remark that, when $q \in S^1$, the existence of such a representation
implies that
the factorials do not vanish, hence that $q$ is not a root of unity.\par
%
%
{}Following Rideau and Winternitz~\cite{rideau}, let us write
$\Psi^J_{MNq}(z,\zb)$
as
\begin{equation}
  \Psi^J_{MNq}(z,\zb) = N^J_{MNq} Q_{Jq}(\eta)\, q^{-NM/2} R^J_{MNq}(\eta)\,
  \zb^{M+N}, \qquad \eta = z \zb.   \label{eq:Psi-z}
\end{equation}
Here, $N^J_{MNq}$ is a constant, which can be expressed as
\begin{eqnarray}
  N^J_{MNq} & = & C_{JNq} \left(\frac{[J+M]_q!}{[J-M]_q!\,
[2J]_q!}\right)^{1/2},
          \nonumber \\
  C_{JNq} & = & \frac{1}{\sqrt{2\pi}} \left(\frac{[J+N]_q!\,
          [2J+1]_q!}{[J-N]_q!}\right)^{1/2} \gamma(J,N,q)   \label{eq:N},
\end{eqnarray}
in terms of some yet undetermined normalization constant $\gamma(J,N,q)$, and
$q$-factorials, defined by $[x]_q! \equiv [x]_q [x-1]_q \ldots [1]_q$ if $x
\in \N^+$,
$[0]_q! \equiv 1$, and $\left([x]_q!\right)^{-1} \equiv 0$ if $x \in \N^-$.
Equation~(\ref{eq:Psi-z}) also contains two functions of~$\eta$, $Q_{Jq}(\eta)$
and $R^J_{MNq}(\eta)$. The latter is a polynomial, whose explicit form is
given by
\begin{eqnarray}
  R^J_{MNq}(\eta) & = & [J-N]_q!\, [J-M]_q! \nonumber \\
  & & \mbox{} \times \sum_k \frac{(-\eta)^k}{[k]_q!\, [J-M-k]_q!\, [J-N-k]_q!\,
          [M+N+k]_q!},   \label{eq:R}
\end{eqnarray}
the summation over $k$ being restricted by the condition that all the
factorials in
the denominator be positive. The former is defined by the functional equation
\begin{equation}
  Q_{Jq}(q^2\eta)(1+\eta) =  Q_{Jq}(\eta) (1+q^{-2J}\eta),      \label{eq:equQ}
\end{equation}
whose solution, only determined up to an arbitrary multiplicative factor
$f_{Jq}(\eta)$ such that
\begin{equation}
  f_{Jq}(q^2\eta) =  f_{Jq}(\eta),      \label{equf}
\end{equation}
will be discussed in detail for both $q \in \R^+$, and generic $q \in S^+$,
in the
next section.
\par
%
%
In terms of spherical coordinates, Eq.~(\ref{eq:Psi-z}) becomes~\cite{rideau}
\begin{eqnarray}
  \Psi^J_{MNq}(\theta,\phi) & = & C_{JNq} \left(\frac{[J-N]_q!}{[J+N]_q!\,
           [2J]_q!}\right)^{1/2} i^{-2J+M+N} q^{-NM/2} \nonumber \\
  & & \mbox{} \times P^J_{MNq}(\cos\theta)\, e^{-i(M+N)\phi},
\label{eq:Psi-theta}
\end{eqnarray}
where
\begin{eqnarray}
  P^J_{MNq}(\xi) & = & i^{2J-M-N} \left(\frac{[J+M]_q!\,
[J+N]_q!}{[J-M]_q!\, [J-N]_q!}
          \right)^{1/2} \eta^{(M+N)/2} Q_{Jq}(\eta)\, R^J_{MNq}(\eta),
\nonumber \\
  & & \xi = \cos\theta, \qquad \eta = \frac{1+\xi}{1-\xi} = \cot^2
\frac{\theta}{2},
          \label{eq:q-Vil}
\end{eqnarray}
are $q$-Vilenkin functions. For integer $J$~values, the functions
$\Psi^J_{M0q}(\theta,\phi)$ are proportional to $q$-spherical harmonics, while
$P_{Jq}(\xi) \equiv P^J_{M0q}(\xi)$ are $q$-analogues of Legendre
polynomials.\par
%
%
In the $q \to 1$ limit, the su$_q$(2) realization~(\ref{eq:su-q}) goes over
into the
su(2) realization
\begin{equation}
  H_3 = - z \partial_z + \zb \partial_{\zb} - N, \qquad H_+ = - \partial_z
- \zb^2
  \partial_{\zb} + N \zb, \qquad H_- = z^2 \partial_z + \partial_{\zb} + N z,
  \label{eq:su}
\end{equation}
the constant $\gamma(J,N,q)$ into $\gamma(J,N,1) = 1$, and the $q$-Vilenkin
functions into ordinary ones $P^J_{MN}(\xi)$. The latter are given by
Eq.~(\ref{eq:q-Vil}), where $[x]_q \to x$, and $Q_{Jq}(\eta) \to Q_J(\eta) =
(1+\eta)^{-J}$. The operators~(\ref{eq:su}) satisfy Eq.~(\ref{eq:hermite}),
and the functions $\Psi^J_{MN}$, $J=|N|$, $|N|+1$,~$\ldots$, $M = -J$, $-J+1$,
$\ldots$,~$J$, form an orthonormal set with respect to the scalar product
\begin{equation}
  \langle \psi_1|\psi_2 \rangle = 2 \int \frac{dz d\zb}{(1+z\zb)^2}\,
  \overline{\psi_1\left(z,\zb\right)}\, \psi_2\left(z,\zb\right) = \frac{1}{2}
  \int_0^{\pi} d\theta\, \sin\theta \int_0^{2\pi} d\phi\,
  \overline{\psi_1\left(\theta,\phi\right)}\, \psi_2\left(\theta,\phi\right),
  \label{eq:prodscal}
\end{equation}
where the integral over $z$,~$\zb$ extends over the whole complex plane.\par
%
%
\section{DETERMINATION OF $Q_{Jq}(\eta)$}
\label{sec:Q}
\setcounter{equation}{0}
{}Following Rideau and Winternitz~\cite{rideau}, as a solution of
Eq.~(\ref{eq:equQ}),
we may consider the function
\begin{equation}
 Q_{Jq}(\eta)  =  {}_1\Phi_0\left(q^{2J}; -; q^2, -q^{-2J}\eta\right) =
          {}_1\Phi_0\left(q^{-2J}; -; q^{-2}, -q^{-2}\eta\right),
\label{eq:Q}
\end{equation}
where ${}_1\Phi_{0}$ is a basic hypergeometric series in the notations of
Ref.~\cite{gasper}.\par
%
%
{}For $q \in \R^+$, use of the $q$-binomial theorem~\cite{gasper} leads to the
expressions
\begin{equation}
  Q_{Jq}(\eta) = \prod_{k=0}^{\infty} \frac{(1+q^{2k}\eta)}{(1+q^{-2J+2k}\eta)}
  \label{eq:Q<}
\end{equation}
if $0<q<1$, and
\begin{equation}
  Q_{Jq}(\eta) = \prod_{k=0}^{\infty} \frac{(1+q^{-2J-2k-2}\eta)}
  {(1+q^{-2k-2}\eta)}   \label{eq:Q>}
\end{equation}
if $q>1$. For integer $J$~values, both expressions reduce to the inverse of a
polynomial,
\begin{equation}
  Q_{Jq}(\eta) = \prod_{k=0}^{J-1}\left( \frac{1}{1+\eta q^{-2J+2k}}\right),
  \label{eq:QJentier}
\end{equation}
whereas for half-integer $J$~values, we are left with convergent infinite
products.\par
%
%
{}For generic $q \in S^1$ and integer $J$~values, Eq.~(\ref{eq:QJentier}) still
remains a valid solution of Eq.~(\ref{eq:equQ}). However, for half-integer
$J$~values, the infinite products contained in Eqs.~(\ref{eq:Q<})
and~(\ref{eq:Q>}), as well as other  expressions of ${}_1\Phi_{0}$ in terms of
infinite series or products, found in Refs.~\cite{exton,gasper}, are
divergent. We
have therefore to look for another solution to Eq.~(\ref{eq:equQ}).\par
%
%
For such a purpose, let us linearize Eq.~(\ref{eq:equQ}) into
\begin{equation}
  K_{Jq}(q^2 \eta) - K_{Jq}(\eta) = \ln \frac{1 + q^{-2J}\eta}{1 + \eta},
  \label{eq:equK}
\end{equation}
by setting
\begin{equation}
  K_{Jq}(\eta) = \ln Q_{Jq}(\eta).
\end{equation}
In terms of the operator $X \equiv \eta \partial_{\eta}$,
Eq.~(\ref{eq:equK}) can be
rewritten as
\begin{equation}
  \left(q^{2X} - 1\right) K_{Jq}(\eta) = \left(q^{-2JX} - 1\right) \ln(1 +
\eta).
  \label{eq:equKbis}
\end{equation}
\par
%
%
Let us consider the difference equation
\begin{equation}
  \left(q^X - q^{-X}\right) L_q(\eta) \equiv L_q(q\eta) - L_q(q^{-1}\eta) =
  \ln(1 + \eta).
  \label{eq:equL}
\end{equation}
If we are able to find a solution to the latter, then
\begin{equation}
  K_{Jq}(\eta) \equiv q^{-X} \left(q^{-2JX} - 1\right) L_q(\eta) =
  L_q\left(q^{-2J-1}\eta\right) - L_q(q^{-1}\eta)
\end{equation}
will be a solution of Eq.~(\ref{eq:equKbis}).\par
%
%
We will now proceed to demonstrate
%
\begin{lemma}    \label{lem-repint}
{}For $0 < \eta < \infty$, and $q = e^{i\tau}$ different from a root of
unity, the
function
\begin{eqnarray}
  L_q(\eta) & = & \frac{1}{2\pi i} \int_0^{\infty} \frac{dt}{t(1+t)} \ln
\left(1 + \eta
          t^{\tau/\pi}\right), \qquad \mbox{if $0 < \tau < \pi$},  \label{eq:Lsol1} \\
  L_q(\eta) & = & - \frac{1}{2\pi i} \int_0^{\infty} \frac{dt}{t(1+t)} \ln
\left(1 +
          \eta t^{-\tau/\pi}\right), \qquad \mbox{if $-\pi < \tau < 0$},
\label{eq:Lsol2}
\end{eqnarray}
is a solution of Eq.~(\ref{eq:equL}).
\end{lemma}
\par
%
%
{\it Proof.} We note that if some function $L_q(\eta)$ is a solution of
Eq.~(\ref{eq:equL}) for $q = e^{i\tau}$, $0 < \tau < \pi$, then $-
L_{q^{-1}}(\eta)$ is
also a solution of the same. Hence, Eq.~(\ref{eq:Lsol2}) directly results from
Eq.~(\ref{eq:Lsol1}). It is also a simple matter to show that the integral
on the
right-hand side of Eq.~(\ref{eq:Lsol1}) is convergent. It therefore only
remains to
prove that the latter satisfies Eq.~(\ref{eq:equL}). For such a purpose, we
have to
separately consider the integral when $\eta$ is replaced by $\eta
e^{i\tau}$, or by
$\eta e^{-i\tau}$.\par
%
%
Let us introduce a function $M(v)$ of a complex variable~$v$, defined by
\begin{equation}
  M(v) = [v(1+v)]^{-1} \ln\left(1 + \eta e^{-i\tau} v^{\tau/\pi}\right),
\label{eq:M}
\end{equation}
where on the right-hand side, there appear two multivalued functions
$v^{\tau/\pi}$, and $\ln w$, where $w = 1 + \eta e^{-i\tau} v^{\tau/\pi}$.\par
%
%
{}For the function $v^{\tau/\pi}$, let us choose a branch cut along the
positive real
axis, so that $v^{\tau/\pi} = |v^{\tau/\pi}| \exp(i\tau\alpha/\pi)$, where
$v = |v|
\exp(i\alpha)$, and $0 < \alpha < 2\pi$. On the two sides of such a cut,
the argument
of the logarithm on the right-hand side of Eq.~(\ref{eq:M}), takes the values
\begin{equation}
  w_+ = 1 + \eta e^{-i\tau} |v|^{\tau/\pi} \qquad \mbox{if $\Re v > 0$,
$\Im v = +0$},
  \label{eq:w+}
\end{equation}
and
\begin{equation}
  w_- = 1 + \eta e^{i\tau} |v|^{\tau/\pi} \qquad \mbox{if $\Re v > 0$, $\Im
v = -0$},
  \label{eq:w-}
\end{equation}
respectively.\par
%
%
Considering next the function $\ln w$, it is easy to show that its branch
point at
$w=0$, and its branch cut along the negative real axis in the complex $w$~plane
cannot be reached within the truncated $v$~plane, since the condition
$\exp\left[i
\tau (\alpha/\pi - 1)\right] = -1$ cannot be fulfilled for $0 < \tau < \pi$,
and $0 < \alpha < 2\pi$.\par
%
%
Hence, when integrating the function $M(v)$ in the complex $v$~plane, one should
consider contours avoiding the branch point $v=0$, the branch cut $\Re
v>0$, $\Im
v=0$, and the simple pole at $v=-1$. Let us consider the two vanishing coutour
integrals
\begin{equation}
  \int_{\Gamma^+} M(v) dv =  \int_{\Gamma^-} M(v) dv = 0,   \label{eq:contour}
\end{equation}
where $\Gamma^+$ and $\Gamma^-$ are the paths in the upper and lower halves of
the $v$~plane, displayed on Fig.~1. The former consists of the upper
half~$C_A^+$
of a large circle of radius~$A$ centred at the origin, and described in the
counterclockwise sense, the upper halves $C_a^+$, $C_{a'}^{'+}$ of two
small circles
of radius $a$,~$a'$, centred at $v=0$ and $v=-1$, respectively, both
described in
the clockwise sense, and three straight lines $L_1^+$, $L_2^+$, $L_3^+$
lying just
above the real axis, and going from $-A$ to $-1-a'$, from $-1+a'$ to $-a$,
and from
$a$ to~$A$, respectively. The latter path~$\Gamma^-$ is defined in a similar
way.\par
%
%
Taking now Eqs.~(\ref{eq:w+}) and (\ref{eq:w-}) into account, we obtain
\begin{equation}
  2\pi i \left(L_q(q\eta) - L_q(q^{-1}\eta)\right) = \lim_{\scriptstyle a
\to 0 \atop
  \scriptstyle A \to \infty} \left\{\int_{L_3^-} M(v) dv - \int_{L_3^+} M(v)
  dv\right\}.   \label{eq:limit}
\end{equation}
Owing to Eq.~(\ref{eq:contour}), each of the integrals on the right-hand side of
Eq.~(\ref{eq:limit}) can be rewritten in terms of integrals along the other
parts of
the path $\Gamma^-$ or~$\Gamma^+$. Those along $L_1^+$ (resp.~$L_2^+$) and
$L_1^-$ (resp.~$L_2^-$) obviously cancel. Furthermore
\begin{equation}
  \lim_{A \to \infty} \left|\int_{C_A^+} M(v) dv - \int_{C_A^-} M(v)
dv\right| \sim
  \lim_{A \to \infty} \frac{\ln A}{A} = 0,
\end{equation}
and
\begin{equation}
  \lim_{a \to 0} \left|\int_{C_a^+} M(v) dv - \int_{C_a^-} M(v) dv\right| \sim
  \lim_{a \to 0} a^{\tau/\pi} = 0,
\end{equation}
so that
\begin{equation}
  2\pi i \left(L_q(q\eta) - L_q(q^{-1}\eta)\right) = - \int_{C'_{a'}} M(v)
  dv = - 2\pi i \Res M(-1) = 2\pi i \ln(1 + \eta),    \label{eq:cqfd}
\end{equation}
where $C'_{a'}$ denotes the circle of radius~$a'$ centred at $v=-1$, and
described in
the counterclockwise sense. Equation~(\ref{eq:cqfd}) completes the
proof.\cqfd\par
%
%
The results of the present section can be collected into
%
\begin{proposition}
The function $Q_{Jq}(\eta)$, appearing on the right-hand side of
Eq.~(\ref{eq:Psi-z}),
is given by Eq.~(\ref{eq:QJentier}) for integer $J$~values, and either $q
\in \R^+$
or generic $q \in S^1$, and by Eqs.~(\ref{eq:Q<}) and~(\ref{eq:Q>}) for
half-integer
$J$~values, and $q \in \R^+$. For half-integer $J$~values, and generic $q
\in S^1$,
it can be expressed as
\begin{equation}
  Q_{Jq}(\eta) = \exp\left\{L_q\left(q^{-2J-1}\eta\right) -
  L_q\left(q^{-1}\eta\right)\right\},  \label{eq:QJdemi-ent}
\end{equation}
where $L_q(\eta)$ admits the integral representation given in
Lemma~\ref{lem-repint}.
\end{proposition}
\par
%
%
\section{UNITARIZATION OF THE REPRESENTATIONS OF su$_q$(2) ON $S^2$}
\label{sec:unitarization}
\setcounter{equation}{0}
In the present section, we will determine a new scalar product $\langle \psi_1|
\psi_2 \rangle_q$ that unitarizes the realization~(\ref{eq:su-q}) of~su$_q$(2),
and goes over into the old one $\langle \psi_1| \psi_2 \rangle$, defined in
Eq.~(\ref{eq:prodscal}), whenever $q \to 1$. For such a purpose, we shall first
impose that Eq.~(\ref{eq:hermite}) is satisfied by the
realization~(\ref{eq:su-q})
with respect to $\langle \psi_1| \psi_2 \rangle_q$. The residual
arbitrariness in
the measure will then be lifted by demanding that  $\langle \psi_1| \psi_2
\rangle_q$ satisfies the usual properties of a scalar product.\par
%
%
We shall successively consider hereunder the cases where $q \in \R^+$, and
generic
$q \in S^1$.\par
%
%
\subsection{The case where $q \in \mbox{\fourteenof R}^+$}
Let us make the following ansatz for $\langle \psi_1| \psi_2 \rangle_q$,
\begin{eqnarray}
  \langle \psi_1| \psi_2 \rangle_q & = & \int_0^{\infty} d\rho
\int_0^{2\pi} d\phi
          \biggl( \overline{A_q \psi_1(\rho,\phi,q)}\, f_1(\rho,q)\, q^{a_1 \rho
          \partial_{\rho}} \psi_2(\rho,\phi,q) \nonumber \\
  & & \mbox{} + \overline{\psi_1(\rho,\phi,q)}\, f_2(\rho,q)\, q^{a_2 \rho
          \partial_{\rho}} A_q \psi_2(\rho,\phi,q)\biggr),    \label{eq:ansatz}
\end{eqnarray}
in terms of the polar coordinates $\rho$,~$\phi$ on the real plane, defined in
Eq.~(\ref{eq:coord}). Here $a_1$,~$a_2$, and $f_1(\rho,q)$, $f_2(\rho,q)$
are some
yet undetermined constants and functions of the indicated arguments,
respectively, and $A_q \equiv q^{-2q \partial_q}$ is the operator that
changes $q$
into~$q^{-1}$, when acting on any function of~$q$,
\begin{equation}
  A_q \psi(\rho,\phi,q) = \psi\left(\rho,\phi,q^{-1}\right).
\end{equation}
\par
%
%
It is easy to check that
\begin{equation}
  \langle \psi_1| H_3 \psi_2 \rangle_q = \langle H_3 \psi_1| \psi_2 \rangle_q
  \label{eq:H_3}
\end{equation}
with respect to~(\ref{eq:ansatz}). Let us now impose the condition
\begin{equation}
  \langle \psi_1| H_+ \psi_2 \rangle_q = \langle H_- \psi_1| \psi_2 \rangle_q.
  \label{eq:H_+-}
\end{equation}
\par
%
%
By combining Eqs.~(\ref{eq:su-q}) and~(\ref{eq:ansatz}), the left-hand side
of this
condition can be written as
\begin{eqnarray}
  \lefteqn{\langle \psi_1| H_+ \psi_2 \rangle_q = \left(q - q^{-1}\right)^{-1}
          \int_0^{\infty} d\rho \int_0^{2\pi} d\phi} \nonumber \\
  & & \mbox{} \times \Biggl\{\overline{\psi_1\left(\rho,\phi,q^{-1}\right)}\,
          f_1(\rho,q)\, e^{-i\phi} \Biggl( - \biggl(q^{a_1}\rho +
          \frac{1}{q^{a_1}\rho}\biggr) q^{\rho\partial_{\rho} - (N/2)}
\nonumber \\
  & & \mbox{} + q^{a_1} \rho\, q^{-i\partial_{\phi} + (3N/2)}
          + \frac{1}{q^{a_1}\rho}\, q^{i\partial_{\phi} - (N/2)} \Biggr)
          \psi_2\left(q^{a_1}\rho,\phi,q\right) \nonumber \\
  & &  - \overline{\psi_1(\rho,\phi,q)}\, f_2(\rho,q)\, e^{-i\phi} \Biggl(
          - \biggl(q^{a_2}\rho + \frac{1}{q^{a_2}\rho}\biggr)
          q^{-\rho\partial_{\rho} + (N/2)} \nonumber \\
  & & \mbox{} + q^{a_2} \rho\, q^{i\partial_{\phi} - (3N/2)} +
         \frac{1}{q^{a_2}\rho}\, q^{-i\partial_{\phi} + (N/2)} \Biggr)
          \psi_2\left(q^{a_2}\rho,\phi,q^{-1}\right) \Biggr\}.
\end{eqnarray}
After integrating by parts and making some straightforward transformations, it
becomes
\begin{eqnarray}
  \lefteqn{\langle \psi_1| H_+ \psi_2 \rangle_q = \left(q - q^{-1}\right)^{-1}
          \int_0^{\infty} d\rho \int_0^{2\pi} d\phi\, e^{-i\phi}} \nonumber \\
  & & \mbox{} \times \Biggl\{ - \Biggl(\biggl(q^{a_1-1}\rho +
          \frac{1}{q^{a_1-1}\rho}\biggr) f_1(q^{-1}\rho,q) q^{-
\rho\partial_{\rho} - 1
          - (N/2)}  \overline{\psi_1\left(\rho,\phi,q^{-1}\right)} \Biggr)
          \psi_2\left(q^{a_1}\rho,\phi,q\right) \nonumber \\
  & & \mbox{} + \Biggl(\biggl(\rho\, q^{i\partial_{\phi} + a_1 + 1 + (3N/2)} +
          \frac{1}{\rho} q^{- i\partial_{\phi} - a_1 - 1 - (N/2)}\biggr)
\overline{\psi_1
          \left(\rho,\phi,q^{-1}\right)} \Biggr) f_1(\rho,q)
          \psi_2\left(q^{a_1}\rho,\phi,q\right) \nonumber \\
  & & \mbox{} + \Biggl(\biggl(q^{a_2+1}\rho + \frac{1}{q^{a_2+1}\rho}\biggr)
          f_2(q\rho,q) q^{\rho\partial_{\rho} + 1 + (N/2)}  \overline{\psi_1(
          \rho,\phi,q)} \Biggr) \psi_2\left(q^{a_2}\rho,\phi,q^{-1}\right)
          \nonumber \\
  & & \mbox{} - \Biggl(\biggl(\rho\, q^{- i\partial_{\phi} + a_2 - 1 - (3N/2)} +
          \frac{1}{\rho} q^{i\partial_{\phi} - a_2 + 1 + (N/2)}\biggr)
\overline{\psi_1
          (\rho,\phi,q)} \Biggr) f_2(\rho,q) \nonumber \\
  & &  \mbox{} \times \psi_2\left(q^{a_2}\rho,\phi,q^{-1}\right)
          \Biggr\}.   \label{eq:H_+}
\end{eqnarray}
\par
%
%
On the other hand, for real $q$~values the right-hand side of
Eq.~(\ref{eq:H_+-})
can be written as
\begin{eqnarray}
  \lefteqn{\langle H_- \psi_1| \psi_2 \rangle_q = \left(q - q^{-1}\right)^{-1}
          \int_0^{\infty} d\rho \int_0^{2\pi} d\phi\, e^{-i\phi}} \nonumber \\
  & & \mbox{} \times \Biggl\{\Biggl(\biggl\{- \biggl(\rho + \frac{1}{\rho}
\biggr)
          q^{- \rho \partial_{\rho} - (N/2)} + \rho\, q^{i \partial_{\phi}
+ (3N/2)} +
          \frac{1}{\rho}\, q^{- i \partial_{\phi} - (N/2)}\biggr\}
          \overline{\psi_1\left(\rho,\phi,q^{-1}\right)}\Biggr) \nonumber \\
  & & \mbox{} \times f_1(\rho,q) \psi_2\left(q^{a_2} \rho, \phi,q\right)
\nonumber
          \\
  & & \mbox{} + \Biggl(\biggl\{\biggl(\rho + \frac{1}{\rho} \biggr) q^{\rho
          \partial_{\rho} + (N/2)} - \rho\, q^{- i \partial_{\phi} - (3N/2)} -
          \frac{1}{\rho}\, q^{i \partial_{\phi} + (N/2)}\biggr\}
          \overline{\psi_1(\rho,\phi,q)}\Biggr) \nonumber \\
  & & \mbox{} \times f_2(\rho,q) \psi_2\left(q^{a_2} \rho, \phi,q^{-1}\right)
          \Biggr\}.    \label{eq:H_-}
\end{eqnarray}
\par
%
%
It now remains to equate the right-hand side of Eq.~(\ref{eq:H_+}) with that of
Eq.~(\ref{eq:H_-}). Both of them being some linear combinations of four
different
types of terms, containing one of the operators $q^{-i\partial_{\phi}}$,
$q^{i\partial_{\phi}}$, $q^{-\rho \partial_{\rho}}$, or~$q^{\rho
\partial_{\rho}}$,
acting on some function, respectively, it is sufficient to separately
equate such
terms. The conditions on the first two classes of terms impose that
\begin{equation}
  a_1 = - 1, \qquad a_2 = 1,   \label{eq:a}
\end{equation}
while those on the last two lead to the equations
\begin{eqnarray}
  q^{-1} \left(q^{-2} \rho + \frac{1}{q^{-2}\rho}\right)
f_1\left(q^{-1}\rho,q\right)
           & = & \left(\rho + \frac{1}{\rho}\right) f_1(\rho,q), \nonumber \\
  q \left(q^2 \rho + \frac{1}{q^2\rho}\right) f_2(q\rho,q)
           & = & \left(\rho + \frac{1}{\rho}\right) f_2(\rho,q),
\end{eqnarray}
whose solutions are given by
\begin{equation}
  f_1(\rho,q) = \frac{B_1(q) q^{-1} \rho}{\left(1+\rho^2\right)
\left(1+q^{-2}\rho^2
  \right)}, \qquad f_2(\rho,q) = \frac{B_2(q) q \rho}{\left(1+\rho^2\right)
  \left(1+q^2\rho^2\right)},     \label{eq:f}
\end{equation}
in terms of two undetermined constants $B_1(q)$, and~$B_2(q)$.\par
%
%
Let us now further restrict the sesquilinear form~(\ref{eq:ansatz}), where
substitutions~(\ref{eq:a}) and~(\ref{eq:f}) have been made, by imposing
that it is
Hermitian, i.e.,
\begin{equation}
  \overline{\langle \psi_1 | \psi_2 \rangle_q} = \langle \psi_2 | \psi_1
\rangle_q.
  \label{eq:hermitebis}
\end{equation}
By a straightforward calculation, similar to that carried out for
condition~(\ref{eq:H_+-}), it can be shown that Eq.~(\ref{eq:hermitebis})
leads to
the relation
\begin{equation}
  B_2(q) = \overline{B_1(q)}.    \label{eq:B}
\end{equation}
As a consequence, there only remains a single undetermined constant $B(q) \equiv
B_1(q)$ in Eq.~(\ref{eq:ansatz}). At this stage, it is important to notice
that had we
only considered a single term, instead of two, in Eq.~(\ref{eq:ansatz}), it
would
have been impossible to fulfil condition~(\ref{eq:hermitebis}).\par
%
%
In addition, we remark that Eqs.~(\ref{eq:H_+-}) and~(\ref{eq:hermitebis}) imply
that
\begin{equation}
  \langle \psi_1 | H_- \psi_2 \rangle_q = \langle H_+ \psi_1 | \psi_2 \rangle_q.
\end{equation}
Hence, all the Hermiticity conditions~(\ref{eq:hermite}) on the su$_q$(2)
generators are satisfied by the form defined in Eqs.~(\ref{eq:ansatz}),
(\ref{eq:a}),
(\ref{eq:f}), and~(\ref{eq:B}). The functions $\Psi^J_{MNq}(z,\zb)$, defined in
Eq.~(\ref{eq:Psi-z}), and corresponding to a fixed $N$~value, but different $J$
and/or $M$~values, are therefore orthogonal with respect to such a form.\par
%
%
To make $\langle \psi_1 | \psi_2 \rangle$ into a scalar product, it only
remains to
impose that it is a positive definite form. Since we also want that in the
resulting
Hilbert space, the functions $\Psi^J_{MNq}$ with given $J$ and $N$~values,
and $M
= -J$, $-J+1$, $\ldots$,~$J$, form an orthonormal basis for the su$_q$(2) irrep
characterized by~$J$, a condition that combines both requirements is
\begin{equation}
  \left\langle \Psi^J_{MNq} | \Psi^J_{MNq} \right\rangle = 1, \qquad M =
-J, -J+1,
  \ldots, J.   \label{eq:norm-M}
\end{equation}
By using Eqs.~(\ref{eq:casimir}) and~(\ref{eq:irrep}) for $M \ne J$,
Eq.~(\ref{eq:norm-M}) can be transformed into the condition
\begin{equation}
  \left\langle \Psi^J_{JNq} | \Psi^J_{JNq} \right\rangle = 1.
\label{eq:norm-J}
\end{equation}
\par
%
%
In Appendix A, the squared norm of~$\Psi^J_{JNq}$ is calculated by using
Eqs.~(\ref{eq:Psi-z}), (\ref{eq:N}), (\ref{eq:R}), (\ref{eq:Q<}), (\ref{eq:Q>}), and by
taking Eqs.~(\ref{eq:ansatz}), (\ref{eq:a}), (\ref{eq:f}), and~(\ref{eq:B})
into account.
The resulting condition~(\ref{eq:norm-J}) reads
\begin{equation}
  \frac{\ln q}{q-q^{-1}} \left(B(q)\, \overline{\gamma(J,N,q^{-1})} \,
\gamma(J,N,q)
  + \overline{B(q)}\, \overline{\gamma(J,N,q)} \,
\gamma(J,N,q^{-1})\right) = 1.
  \label{eq:norm-J-bis}
\end{equation}
Since in the limit $q \to 1$, $\gamma(J,N,q) \to 1$, we may choose
\begin{equation}
  \gamma(J,N,q) = 1, \qquad B(q) = \overline{B(q)} = \frac{q-q^{-1}}{2\ln q} =
  \frac{\sinh\tau}{\tau}.
\end{equation}
{}For $q\to 1$ or $\tau\to 0$, we find that $B(q) \to 1$, so that $\langle
\psi_1 |
\psi_2 \rangle_q \to \langle \psi_1 | \psi_2 \rangle$, where the latter is given
by Eq.~(\ref{eq:prodscal}), as it should be.\par
%
%
The results obtained can be summarized as follows:
%
\begin{proposition}   \label{prop-prodscal}
{}For $q \in \R^+$, the scalar product
\begin{eqnarray}
  \langle \psi_1 | \psi_2 \rangle_q & = & \frac{q-q^{-1}}{2\ln q} \int dz\, d\zb
         \Biggl( \overline{\psi_1(z,\zb,q^{-1})}\, \frac{1}{(1+z\zb)
(1+q^{-2}z\zb)}\,
         q^{- z\partial_z - \zb\partial_{\zb} - 1} \psi_2(z,\zb,q) \nonumber \\
  & & \mbox{} + \overline{\psi_1(z,\zb,q)}\, \frac{1}{(1+z\zb) (1+q^2 z\zb)}\,
         q^{z\partial_z + \zb\partial_{\zb} + 1} \psi_2(z,\zb,q^{-1}) \Biggr),
\end{eqnarray}
or
\begin{eqnarray}
  \lefteqn{\langle \psi_1 | \psi_2 \rangle_q = \frac{q-q^{-1}}{8\ln q}
\int_0^{\pi}
         d\theta\, \sin\theta \int_0^{2\pi} d\phi} \nonumber \\
  & & \mbox{} \times \Biggl(\overline{\psi_1(\theta,\phi,q^{-1})}\, \frac{1}
         {\sin^2(\theta/2) + q^{-2} \cos^2(\theta/2)}\,
q^{\sin\theta\partial_{\theta}
         - 1} \psi_2(\theta,\phi,q) \nonumber \\
  & & \mbox{} + \overline{\psi_1(\theta,\phi,q)}\, \frac{1}{\sin^2(\theta/2)
         + q^2 \cos^2(\theta/2)}\, q^{- \sin\theta\partial_{\theta} + 1}
         \psi_2(\theta,\phi,q^{-1}) \Biggr),
\end{eqnarray}
unitarizes the su$_q$(2) realization~(\ref{eq:su-q}), where $N$ may take any
integer or half-integer value. The functions $\Psi^J_{MNq}(z,\zb)$, or
$\Psi^J_{MNq}(\theta,\phi)$, defined in Eqs.~(\ref{eq:Psi-z})
and~(\ref{eq:Psi-theta}), where $J = |N|$, $|N|+1$,~$\ldots$, $M = -J$, $-J+1$,
$\ldots$,~$J$, and $\gamma(J,N,q) = 1$, form an orthonormal set with respect to
such a scalar product.
\end{proposition}
\par
%
%
{}From Proposition~\ref{prop-prodscal}, we easily obtain
%
\begin{corollary}
{}For $q \in \R^+$, the $q$-Vilenkin functions $P^J_{MNq}(\xi)$, defined in
Eq.~(\ref{eq:q-Vil}), satisfy the orthonormality relation
\begin{eqnarray}
  \lefteqn{\frac{q-q^{-1}}{4\ln q} \int_{-1}^{+1} d\xi \Biggl(
         \overline{P^{J'}_{MNq^{-1}}(\xi)}\, \frac{1}{q+q^{-1} -
(q-q^{-1})\xi}\,
         q^{\left(\xi^2-1\right) \partial_{\xi}}\, P^J_{MNq}(\xi)} \nonumber \\
  & & \mbox{} + \overline{P^{J'}_{MNq}(\xi)}\, \frac{1}{q+q^{-1} +
(q-q^{-1})\xi}\,
         q^{-\left(\xi^2-1\right) \partial_{\xi}}\, P^J_{MNq^{-1}}(\xi) \Biggr)
         = \frac{\delta_{J',J}}{[2J+1]_q}.
\end{eqnarray}
\end{corollary}
%
\subsection{The case where $q \in S^1$}
Whenever $q\in S^1$, the ansatz~(\ref{eq:ansatz}) does not work, because though
Eq.~(\ref{eq:H_+}) remains valid, Eq.~(\ref{eq:H_-}) is changed in such a
way that
both cannot be matched. Let us therefore change Eq.~(\ref{eq:ansatz}) into the
following ansatz
\begin{eqnarray}
  \langle \psi_1| \psi_2 \rangle_q & = & \int_0^{\infty} d\rho
\int_0^{2\pi} d\phi
          \biggl( \overline{\psi_1(\rho,\phi,q)}\, f_1(\rho,q)\, q^{a_1 \rho
          \partial_{\rho}} \psi_2(\rho,\phi,q) \nonumber \\
  & & \mbox{} + \overline{A_q \psi_1(\rho,\phi,q)}\, f_2(\rho,q)\, q^{a_2 \rho
          \partial_{\rho}} A_q \psi_2(\rho,\phi,q)\biggr),
\label{eq:ansatzbis}
\end{eqnarray}
where $a_1$, $a_2$, $f_1(\rho,q)$, $f_2(\rho,q)$, and~$A_q$ keep the same
meaning as before.\par
Condition~(\ref{eq:H_3}) is again automatically satisfied. Turning now to
condition~(\ref{eq:H_+-}), it is easy to see that Eqs.~(\ref{eq:H_+})
and~(\ref{eq:H_-}) remain valid, except for the interchange of
$\overline{\psi_1(\rho,\phi,q)}$ with
$\overline{\psi_1\left(\rho,\phi,q^{-1}\right)}$. Hence,
Eq.~(\ref{eq:H_+-}) is also
fulfilled by choosing $a_1$, $a_2$, $f_1(\rho,q)$, and~$f_2(\rho,q)$ as given in
Eqs.~(\ref{eq:a}), and~(\ref{eq:f}).\par
A difference with the case where $q \in \R^+$ appears when imposing the
Hermiticity condition~(\ref{eq:hermitebis}). The latter is now equivalent to the
relations
\begin{equation}
  \overline{B_1(q)} = B_1(q), \qquad \overline{B_2(q)} = B_2(q),
\label{eq:Bbis}
\end{equation}
showing that the real constants $B_1(q)$, and~$B_2(q)$ remain independent.
In the
present case, keeping only one of the two terms on the right-hand side of
Eq.~(\ref{eq:ansatzbis}) would therefore lead to a well-behaved scalar
product.\par
As shown in Appendix B, condition~(\ref{eq:norm-J}) now reads
\begin{equation}
  \frac{\ln q}{q-q^{-1}} \left(B_1(q) |\gamma(J,N,q)|^2 + B_2(q)
  |\gamma(J,N,q^{-1})|^2 \right) = 1. \label{eq:norme S1}
\end{equation}
Among the infinitely many solutions of this equation, we may select the most
symmetrical one,
\begin{equation}
  \gamma(J,N,q) = 1, \qquad B_1(q) = B_2(q) = \frac{q-q^{-1}}{2\ln q} =
  \frac{\sin\tau}{\tau}.
\end{equation}
Hence, whenever $q\to 1$ or $\tau\to 0$, the limit of $\langle \psi_1 |\psi_2
\rangle_q$ is again $\langle \psi_1 | \psi_2 \rangle$, as it should be.\par
%
%
In conclusion, we obtain
%
\begin{proposition}  \label{prop-prodscalbis}
{}For generic $q \in S^1$, the scalar product
\begin{eqnarray}
  \langle \psi_1 | \psi_2 \rangle_q & = & \frac{q-q^{-1}}{2\ln q} \int dz\, d\zb
         \Biggl( \overline{\psi_1(z,\zb,q)}\, \frac{1}{(1+z\zb)
(1+q^{-2}z\zb)}\,
         q^{- z\partial_z - \zb\partial_{\zb} - 1} \psi_2(z,\zb,q) \nonumber \\
  & & \mbox{} + \overline{\psi_1(z,\zb,q^{-1})}\, \frac{1}{(1+z\zb) (1+q^2
z\zb)}\,
         q^{z\partial_z + \zb\partial_{\zb} + 1} \psi_2(z,\zb,q^{-1}) \Biggr),
\end{eqnarray}
or
\begin{eqnarray}
  \lefteqn{\langle \psi_1 | \psi_2 \rangle_q = \frac{q-q^{-1}}{8\ln q}
\int_0^{\pi}
         d\theta\, \sin\theta \int_0^{2\pi} d\phi} \nonumber \\
  & & \mbox{} \times \Biggl(\overline{\psi_1(\theta,\phi,q)}\, \frac{1}
         {\sin^2(\theta/2) + q^{-2} \cos^2(\theta/2)}\,
q^{\sin\theta\partial_{\theta}
         - 1} \psi_2(\theta,\phi,q) \nonumber \\
  & & \mbox{} + \overline{\psi_1(\theta,\phi,q^{-1})}\,
\frac{1}{\sin^2(\theta/2)
         + q^2 \cos^2(\theta/2)}\, q^{- \sin\theta\partial_{\theta} + 1}
         \psi_2(\theta,\phi,q^{-1}) \Biggr),
\end{eqnarray}
unitarizes the su$_q$(2) realization~(\ref{eq:su-q}), where $N$ may take any
integer or half-integer value. The functions $\Psi^J_{MNq}(z,\zb)$, or
$\Psi^J_{MNq}(\theta,\phi)$, defined in Eqs.~(\ref{eq:Psi-z})
and~(\ref{eq:Psi-theta}), where $J = |N|$,
$|N|+1$,~$\ldots$, $M = -J$, $-J+1$, $\ldots$,~$J$, and $\gamma(J,N,q) =
1$, form
an orthonormal set with respect to such a scalar product.
\end{proposition}
%
\begin{corollary}    \label{prop-orthobis}
{}For generic $q \in S^1$, the $q$-Vilenkin functions $P^J_{MNq}(\xi)$,
defined in
Eq.~(\ref{eq:q-Vil}), satisfy the orthonormality relation
\begin{eqnarray}
  \lefteqn{\frac{q-q^{-1}}{4\ln q} \int_{-1}^{+1} d\xi \Biggl(
         \overline{P^{J'}_{MNq}(\xi)}\, \frac{1}{q+q^{-1} - (q-q^{-1})\xi}\,
         q^{\left(\xi^2-1\right) \partial_{\xi}}\, P^J_{MNq}(\xi)} \nonumber \\
  & & \mbox{} + \overline{P^{J'}_{MNq^{-1}}(\xi)}\, \frac{1}{q+q^{-1} +
(q-q^{-1})\xi}\,
         q^{-\left(\xi^2-1\right) \partial_{\xi}}\, P^J_{MNq^{-1}}(\xi) \Biggr)
         = \frac{\delta_{J',J}}{[2J+1]_q}.
\end{eqnarray}
\end{corollary}
\par
%
%
\section{Conclusion}   \label{sec:conclusion}
In the present paper, we did extend the study of the su$_q$(2)
representations on a
real two-dimensional sphere, carried out by Rideau and
Winternitz~\cite{rideau}, in
two ways.\par
%
%
{}Firstly, we did prove that such representations exist not only for $q \in
\R^+$, but
also for generic $q \in S^1$. For such a purpose, we did provide an integral
representation for the functions $Q_{Jq}(\eta)$, entering the definition of the
$q$-Vilenkin functions, whenever $J$ takes any half-integer value.\par
%
%
Secondly, we did unitarize the representations by determining appropriate scalar
products for both ranges of $q$~values. Such scalar products are expressed
in terms
of ordinary integrals, instead of $q$-integrals, as is usually the
case~\cite{vilenkin2}.\par
%
%
The resulting orthonormality relations for the $q$-Vilenkin and related
functions
should play an important role in applications to quantum mechanics, such as
those
considered in Refs.~\cite{irac1,irac2}.\par
%
%
\newpage
\section*{APPENDIX A: PROOF OF EQUATION~(\ref{eq:norm-J-bis})}
\label{sec:appendixA}
\renewcommand{\theequation}{A\arabic{equation}}
\setcounter{equation}{0}
The purpose of this appendix is to evaluate the squared norm of the function
$\Psi^J_{JNq}(z,\zb)$ when the scalar product~(\ref{eq:ansatz}) is used, and
Eqs.~(\ref{eq:a}), (\ref{eq:f}), and~(\ref{eq:B}) are taken into account.\par
%
%
{}From Eqs.~(\ref{eq:complex}), (\ref{eq:Psi-z}), (\ref{eq:N}),
and~(\ref{eq:R}),
$\Psi^J_{JNq}$ can be written in polar coordinates as
\begin{equation}
  \Psi^J_{JNq} = \frac{C_{JNq}}{[J+N]_q!}\, q^{-JN/2} Q_{Jq}(\rho^2)\,
\rho^{J+N}
  e^{-i(J+N)\phi}.  \label{eq:Psi-M=J}
\end{equation}
Its squared norm can therefore be expressed as
\begin{eqnarray}
  \Bigl\langle \Psi^J_{JNq} \Big| \Psi^J_{JNq} \Bigr\rangle_q & = & \frac{\pi}
          {([J+N]_q!)^2} \biggl(B(q)\, \overline{C_{JNq^{-1}}}\, C_{JNq}\,
q^{-J-N-1}
          {\cal I}_q \nonumber \\
  & & \mbox{} + \overline{B(q)}\, \overline{C_{JNq}}\, C_{JNq^{-1}}\, q^{J+N+1}
          {\cal I}_{q^{-1}} \biggr),    \label{eq:norm-def}
\end{eqnarray}
in terms of the integral
\begin{equation}
  {\cal I}_q = \int_0^{\infty} d\eta\, Q_{Jq^{-1}}(\eta)
\frac{\eta^{J+N}}{(1+\eta)
  (1+q^{-2}\eta)}\, Q_{Jq}\left(q^{-2} \eta\right),   \label{eq:I-def}
\end{equation}
and the same with $q$ replaced by~$q^{-1}$.\par
%
%
By introducing Eqs.~(\ref{eq:Q<}) and~(\ref{eq:Q>}) into
Eq.~(\ref{eq:I-def}), we
obtain
\begin{equation}
  {\cal I}_q = \int_0^{\infty} d\eta\, \eta^{J+N} \prod_{k=0}^{\infty}
  \frac{(1+q^{2J+2k+2}\eta)}{(1+q^{-2J+2k-2}\eta)} = q^{2(J+1)(J+N+1)}
  \tilde{B}_{q^2}(J+N+1,J-N+1)   \label{eq:I-B<}
\end{equation}
if $0 < q < 1$, and
\begin{equation}
  {\cal I}_q = \int_0^{\infty} d\eta\, \eta^{J+N} \prod_{k=0}^{\infty}
  \frac{(1+q^{-2J-2k-4}\eta)}{(1+q^{2J-2k}\eta)} = q^{-2J(J+N+1)}
  \tilde{B}_{q^{-2}}(J+N+1,J-N+1)   \label{eq:I-B>}
\end{equation}
if $q > 1$. In Eqs.~(\ref{eq:I-B<}) and~(\ref{eq:I-B>}), we denote by
$\tilde{B}_q(x,y)$ Ramanujan's continuous $q$-analogue of the beta
integral~\cite{andrews}
\begin{equation}
  \tilde{B}_q(x,y) = \int_0^{\infty} dt\, t^{x-1} \prod_{k=0}^{\infty}
  \frac{(1+q^{x+y+k}t)}{(1+q^kt)}, \qquad 0<q<1,   \label{eq:betaint}
\end{equation}
to distinguish it from the discrete $q$-analogue of the same, known as
$B_q(x,y)$
(see e.g. Eq.~(1.11.7) of Ref.~\cite{gasper}).\par
%
%
{}From Eq.~(5.8) of Ref.~\cite{andrews}, $\tilde{B}_q(x,y)$ is given for generic
$x$~values by
\begin{equation}
  \tilde{B}_q(x,y) = \frac{\pi}{\sin\pi x} \prod_{k=1}^{\infty}
\frac{(1-q^{k-x})
  (1-q^{x+y+k-1})}{(1-q^k)(1-q^{y+k-1})}.   \label{eq:B_q}
\end{equation}
The values of~$x$, which appear in Eqs.~(\ref{eq:I-B<})
and~(\ref{eq:I-B>}), being
$x=J+N+1 \in \N^+$, we have to calculate the limit of the right-hand side of
Eq.~(\ref{eq:B_q}) when $x \to m \in \N^+$. Using L'Hospital rule, we find
\begin{equation}
  \lim_{x\to m} \frac{1-q^{m-x}}{\sin\pi x} = (-1)^m \frac{\ln q}{\pi},
\qquad m \in
  \N^+.
\end{equation}
Hence, for $x=m$, $y=n$, $m$, $n \in \N^+$, Eq.~(\ref{eq:B_q}) becomes
\begin{eqnarray}
  \tilde{B}_q(m,n) & =  & (-1)^m (\ln q) \frac{\prod_{k=1}^{m-1} (1-q^{k-m})}
          {\prod_{k=1}^m (1-q^{n+k-1})} \nonumber \\
  & = & \frac{(\ln q) q^{-m(n+m-1)/2} [m-1]_{q^{1/2}}!\, [n-1]_{q^{1/2}}!}
          {(q^{1/2}-q^{-1/2}) [n+m-1]_{q^{1/2}}!},   \label{eq:B_qbis}
\end{eqnarray}
where in the last step, we introduced $q$-factorials, defined as in
Sec.~\ref{sec:representations}.\par
%
%
{}From Eqs.~(\ref{eq:I-B<}), (\ref{eq:I-B>}), and~(\ref{eq:B_qbis}), it
follows that for
any $q \in \R^+$
\begin{equation}
  {\cal I}_q = \frac{2 (\ln q) q^{J+N+1} [J+N]_q!\, [J-N]_q!} {(q-q^{-1})
[2J+1]_q!}.
\label{calI}
\end{equation}
By taking Eq.~(\ref{eq:N}) into account, the squared norm of~$\Psi^J_{JNq}$,
defined in Eq.~(\ref{eq:norm-def}), therefore becomes
\begin{equation}
  \Bigl\langle \Psi^J_{JNq} \Big| \Psi^J_{JNq} \Bigr\rangle_q =
  \frac{\ln q}{q-q^{-1}} \left(B(q)\, \overline{\gamma(J,N,q^{-1})} \,
\gamma(J,N,q)
  + \overline{B(q)}\, \overline{\gamma(J,N,q)} \,  \gamma(J,N,q^{-1})\right),
\end{equation}
which proves Eq.~(\ref{eq:norm-J-bis}).\par
%
%
\newpage
\section*{APPENDIX B: PROOF OF EQUATION~(\ref{eq:norme S1})}
\label{sec:appendixB}
\renewcommand{\theequation}{B\arabic{equation}}
\setcounter{equation}{0}
The purpose of this appendix is to evaluate the squared norm of the function
$\Psi^J_{JNq}(z,\zb)$ when the scalar product~(\ref{eq:ansatzbis}) is used, and
Eqs.~(\ref{eq:a}), (\ref{eq:f}), and~(\ref{eq:Bbis}) are taken into account.\par
%
%
Since for $q \in S^1$, $\Psi^J_{JNq}$ is still given by
Eq.~(\ref{eq:Psi-M=J}), its
squared norm reads
\begin{eqnarray}
  \Bigl\langle \Psi^J_{JNq} \Big| \Psi^J_{JNq} \Bigr\rangle_q & = & \frac{\pi}
          {([J+N]_q!)^2} \biggl(B_1(q)\, \Big| C_{JNq}\Big|^2\, q^{-J-N-1}
{\cal I}'_q
          \nonumber \\
  & & \mbox{} + B_2(q) \,\Big|C_{JNq^{-1}}\Big|^2\, q^{J+N+1} {\cal I}'_{q^{-1}}
          \biggr).    \label{eq:norm-defS1}
\end{eqnarray}
Here ${\cal I}'_q$ denotes the integral
\begin{equation}
  {\cal I}'_q = \int_0^{\infty} d\eta\, F_{Jq}(\eta)\, \eta^{J+N}
     \label{eq:I-defS1}
\end{equation}
with
\begin{equation}
  F_{Jq}(\eta) =  \overline { Q_{Jq}(\eta)}\, \frac{1}{(1+\eta)
(1+q^{-2}\eta)}\,
  Q_{Jq}\left(q^{-2} \eta\right).   \label{eq:defF}
\end{equation}
\par
%
%
According to whether $J$ is integer or half-integer, we have to insert
Eq.~(\ref{eq:QJentier}) or Eq.~(\ref{eq:QJdemi-ent}) into
Eq.~(\ref{eq:defF}). In both
cases, the result reads
\begin{equation}
  F_{Jq}(\eta)=\prod _{p=0}^{2J+1}\frac{1}{1+q^{2J-2p}\eta}.    \label{eq:F0}
\end{equation}
This is obvious in the former case. In the latter, by using the property
$\overline{L_q(\eta)} = - L_q(\eta)$, Eq.~(\ref{eq:defF}) can be
transformed into
\begin{eqnarray}
  F_{Jq}(\eta) & = & \exp\left\{- L_q\left(q^{2J+1}\eta\right) +
L_q(q\eta)\right\}
           \frac{1}{(1+\eta)(1+q^{-2}\eta)}
\exp\left\{L_q\left(q^{-2J-3}\eta\right) -
           L_q(q^{-3}\eta)\right\} \nonumber \\
  & = & \frac{1}{(1+\eta)(1+q^{-2}\eta)} \exp\Biggl\{ - \sum_{p=0}^{2J+1} \left[
           L_q\left(q^{2J+1-2p}\eta\right) - L_q(q^{2J-1-2p}\eta)\right]
\nonumber \\
  & & \mbox{} + \left[L_q\left(q\eta\right) - L_q(q^{-1}\eta)\right] +
           \left[L_q\left(q^{-1}\eta\right) - L_q(q^{-3}\eta)\right]\Biggr\}.
 \end{eqnarray}
Repeated use of Eq.~(\ref{eq:equL}) for various arguments then directly
leads to the
searched for result~(\ref{eq:F0}).\par
%
%
To evaluate ${\cal I}'_q$ for $F_{Jq}(\eta)$ given by Eq.~(\ref{eq:F0}), we
cannot
use the same method as that employed in Appendix~A to calculate ${\cal I}_q$,
because in the $q$-analogue of the beta integral, given in
Eq.~(\ref{eq:betaint}), $q$
is assumed real. Let us therefore rewrite the integrand of~${\cal I}'_q$ in
the form
\begin{equation}
  \eta^{J+N} F_{Jq}(\eta) = \sum_{p=0}^{2J+1} \frac{a_p^{(J)}}{\eta +
q^{2p-2J}},
\end{equation}
where the coefficient $a_p^{(J)}$ is the residue of $\eta^{J+N}
F_{Jq}(\eta)$ at the
pole $\eta = -q^{2p-2N}$, i.e.,
\begin{equation}
  a_p^{(J)} = (-1)^{J+N}  \frac{q^{J+1}} {(q-q^{-1})^{2J+1}} \times \frac{(-1)^p
  q^{N(2p-2J)}}{[p]_q!\, [2J-p+1]_q!}.
\end{equation}
Then
\begin{eqnarray}
  G_{Jq}(\eta) & \equiv & \int d\eta\, F_{Jq}(\eta) \eta^{J+N} \nonumber \\
  & = & (-1)^{J+N} \frac{q^{J+1}}{(q-q^{-1})^{2J+1}} \sum_{p=0}^{2J+1}
\frac{(-1)^p
         q^{N(2p-2J)}}{[p]_q!\, [2J-p+1]_q!} \ln(\eta + q^{2p-2J}).
\label{eq:G}
\end{eqnarray}
\par
%
%
To calculate the values of $G_{Jq}(\eta)$ for $\eta \to \infty$ and
$\eta=0$, the
following identities~\cite{exton} will be useful:
\begin{eqnarray}
  (1;\eta)^{2J+1}_q & \equiv & \sum_{p=0}^{2J+1} \bin{2J+1}{p} \eta^p =
          \prod_{p=0}^{2J} \left(1 + q^{2p-2J} \eta\right),
\label{eq:1+eta}\\[0.1cm]
  \frac{d}{d\eta} (1;\eta)^{2J+1}_q & = & \sum_{p=0}^{2J+1} \bin{2J+1}{p} p
          \eta^{p-1} = \sum_{p=0}^{2J} q^{2p-2J} \prod_{\scriptstyle r=0 \atop
          \scriptstyle r\ne p}^{2J} \left(1 + q^{2r-2J} \eta\right),
\label{eq:derive}
\end{eqnarray}
where
\begin{equation}
  \bin{n}{p} \equiv \frac{[n]_q!}{[p]_q!\, [n-p]_q!}
\end{equation}
is a $q$-binomial coefficient. From Eq.~(\ref{eq:1+eta}), we obtain
\begin{equation}
  \left(1;-q^{2N}\right)^{2J+1}_q = \sum_{p=0}^{2J+1} (-1)^p \bin{2J+1}{p}
q^{2Np}
  = 0,   \label{eq:1+etabis}
\end{equation}
because on the right-hand side, the factor $(1 - q^{2p-2J+2N})$ vanishes
for $p =
J-N$. Similarly, from Eq.~(\ref{eq:derive}), we get
\begin{eqnarray}
  \left. \frac{d}{d\eta} (1;\eta)^{2J+1}_q \right|_{\eta=- q^{2N}} & = &
          \sum_{p=0}^{2J+1} (-1)^{p-1} p \bin{2J+1}{p} q^{2N(p-1)}
\nonumber \\[0.1cm]
  & = & (-1)^{J+N} (q - q^{-1})^{2J} [J-N]_q!\, [J+N]_q!\, q^{N(2J-1)},
          \label{eq:derivebis}
\end{eqnarray}
since on the right-hand side, only the term corresponding to $p=J-N$ leads to a
nonvanishing result.\par
%
%
By noting that for $\eta \gg 1$,
\begin{equation}
  \ln(\eta + q^{2p-2J}) \simeq \ln (\eta) + \frac{q^{2p-2J}}{\eta} +
   O\left(\frac{1}{\eta^2}\right),
\end{equation}
it directly results from Eq.~(\ref{eq:1+etabis}) that
\begin{equation}
  \lim_{\eta\to\infty} G_{Jq}(\eta) = 0.    \label{eq:Ginfty}
\end{equation}
Furthermore, from Eqs.~(\ref{eq:1+etabis}) and~(\ref{eq:derivebis}), we obtain
\begin{eqnarray}
  G_{Jq}(0)  & = & \frac{(-1)^{J+N} q^{J+1} \ln q}{(q-q^{-1})^{2J+1}}
   \sum_{p=0}^{2J+1}
          \frac{(-1)^p (2p-2J)q^{N(2p-2J)}}{[p]_q!\, [2J-p+1]_q!} \nonumber
\\[0.1cm]
  & = & \frac{2 (-1)^{J+N} q^{J+1} \ln q}{[2J+1]_q!\, (q-q^{-1})^{2J+1}}
   \sum_{p=0}^{2J+1}
          (-1)^p p \bin{2J+1}{p} q^{N(2p-2J)} \nonumber \\[0.1cm]
  & = & - \frac{2 [J+N]_q![J-N]_q! q^{J+N+1} \ln q}{[2J+1]_q!\, (q - q^{-1})}.
   \label{eq:G0}
\end{eqnarray}
\par
%
%
By taking Eqs.~(\ref{eq:G}), (\ref{eq:Ginfty}), and~(\ref{eq:G0}) into
account, we
conclude that for generic $q \in S^1$, ${\cal I}'_q$, defined in
Eq.~(\ref{eq:I-defS1}),
is given by
\begin{equation}
  {\cal I'}_q = \frac{2 (\ln q) q^{J+N+1} [J+N]_q!\, [J-N]_q!} {(q-q^{-1})
[2J+1]_q!}.
\end{equation}
By combining this result with Eqs.~(\ref{eq:N}), (\ref{eq:norm-J}),
and~(\ref{eq:norm-defS1}), Eq.~(\ref{eq:norme S1}) directly follows.\par
%
%

%
%
\newpage
\parindent0cm
\section*{Figure caption}
{\bf Fig.~1.} Contours in the complex $v$~plane used in the proof of
Lemma~\ref{lem-repint}.


\newpage
\begin{thebibliography}{99}
%
\bibitem{exton} H.~Exton {\em $q$-Hypergeometric Functions and Applications}
(Ellis Horwood, New York, 1983).
%
\bibitem{andrews} G.~E.~Andrews, {\em $q$-Series: Their Development and
Application in Analysis, Number Theory, Combinatorics, Physics, and Computer
Algebra} (American Mathematical Society, Providence, RI, 1986).
%
\bibitem{gasper} G.~Gasper and M.~Rahman, {\em Basic Hypergeometric Series}
(Cambridge University, Cambridge, 1990).
%
\bibitem{chari} V.~Chari and A.~Pressley, {\em A Guide to Quantum Groups}
(Cambridge University, Cambridge, 1994).
%
\bibitem{vilenkin2} N.~Ja.~Vilenkin and A.~U.~Klimyk, {\em Representation of Lie
Groups and Special Functions} (Kluwer, Dordrecht, 1991), 3 Vol.
%
\bibitem{rideau} G.~Rideau and P.~Winternitz, J. Math. Phys. {\bf 34}, 6030
(1993).
%
\bibitem{vilenkin1} N.~Ja.~Vilenkin, {\em Special Functions and the Theory of
Group Representations} (American Mathematical Society, Providence, RI, 1968).
%
\bibitem{footnote} Note that the functions $P^J_{MN}(\cos\theta)$, defined in
Ref.~\cite{vilenkin2}, Vol.~1, p.~85, differ by a phase factor from those
considered
in Refs.~\cite{rideau, vilenkin1}.
%
\bibitem{schmidt} J.~R.~Schmidt, J. Math. Phys. {\bf 37}, 3062 (1996).
%
\bibitem{irac1} M.~Irac-Astaud, Lett. Math. Phys. {\bf 36}, 169 (1996);
Czech. J.
Phys. {\bf 46}, 179 (1996).
%
\bibitem{irac2} M.~Irac-Astaud and C.~Quesne, unpublished.

\end{thebibliography}
\end{document}